\documentclass[a4paper,fleqn]{cas-sc}

\usepackage[authoryear]{natbib}
\usepackage{graphicx,subfigure}
\usepackage{adjustbox} 
\usepackage{pifont} 
\usepackage{lineno} 
\usepackage{lipsum} 

\usepackage{booktabs}  
\usepackage{hyperref}
\usepackage{url}
\def\tsc#1{\csdef{#1}{\textsc{\lowercase{#1}}\xspace}}
\tsc{WGM}
\tsc{QE}

\begin{document}
\let\WriteBookmarks\relax
\def\floatpagepagefraction{1}
\def\textpagefraction{.001}

\shorttitle{PCM-NN for Soybean Pod Borer Dynamics Prediction}    

\shortauthors{X. Chen, W. Li, S. Liu and Y. Gao}  

\title [mode = title]{Modeling the Temperature-Humidity Coupling Dynamics of Soybean Pod Borer Population and Assessing the Predictive Performance of the PCM-NN Algorithm}


%

\author[1]{Xu Chen}
\fnmark[1,$\star$]
\tnotetext[1]{These authors contributed equally to this work.}
\credit{Writing – original draft, Visualization, Validation,
Software, Methodology, Investigation, Formal analysis, Data curation}
\affiliation[1]{organization={School of Artificial Intelligence},
            addressline={Jilin University}, 
            city={Changchun},
            postcode={130012}, 
            state={Jilin},
            country={China}}
\fntext[1]{http://orcid.org/0000-0002-0866-6857}
\author[2]{Wenxuan Li}%
\fnmark[1,$\star$]
\credit{Writing – original draft, Visualization, Validation,
Software, Methodology, Investigation, Formal analysis, Data curation}
\affiliation[2]{organization={School of Mathematics},
            addressline={Jilin University}, 
            city={Changchun},
            postcode={130012}, 
            state={Jilin},
            country={China}}
\fntext[2]{http://orcid.org/0000-0002-3601-117X}
\author[3]{Xiaoshuang Li}%
\fnmark[2]
\credit{Investigation, Formal analysis, Data curation}
\affiliation[3]{organization={College of Plant Protection},
            addressline={Jilin Agricultural University, Changchun}, 
            city={Changchun},
            postcode={130118}, 
            state={Jilin},
            country={China}}
\fntext[3]{http://orcid.org/0009-0002-3855-0289}
\author[2]{Suli Liu}
\fnmark[3]
\cormark[1]
\ead{liusuli@jlu.edu.cn}
\credit{Writing – review \& editing, Writing –
original draft, Supervision, Resources, Project administration, Funding
acquisition, Conceptualization}
\cortext[1]{Corresponding author}
\fntext[4]{http://orcid.org/0000-0002-4354-2665}
\author[4]{Yu Gao}
\fnmark[4]
\cormark[1]
\ead{gaothrips@jlau.edu.cn}
\credit{Writing – review \& editing, Writing –
original draft, Supervision, Resources, Project administration, Funding
acquisition, Conceptualization.}
\affiliation[4]{organization={Key Laboratory of Soybean Disease and Pest Correlation},
            addressline={Ministry of Agriculture and Rural Affairs}, 
            city={Changchun},
            postcode={130118}, 
            state={Jilin},
            country={China}}
\fntext[5]{http://orcid.org/0000-0002-3369-8578}

\begin{abstract}
Against the backdrop of global climate change and agricultural globalization, soybean production is increasingly threatened by pest outbreaks, with Leguminivora glycinivorella (commonly known as the soybean pod borer) being a major pest species. This pest is widely distributed, particularly in northeastern China—the country’s primary soybean-producing region—where its outbreaks have significantly affected both yield and quality. Although statistical and mechanistic models have been applied to pest forecasting, existing approaches often fail to effectively integrate climatic factors with pest dynamics and lack sufficient expressive power. To address these limitations, this study proposes a novel pest prediction method based on Physics-Informed Neural Networks (PINNs). Specifically, we formulate a logistic-type ordinary differential equation (ODE) that incorporates microclimate factors—temperature, humidity, and time—to describe the temporal dynamics of the soybean pod borer population. This ODE model is embedded into the PINN framework to develop the Pest Correlation Model Neural Network (PCM-NN), which is used to jointly infer the microclimate-driven parameter function $\alpha(T, H, t)$ and fit the pest population dynamics. We evaluate PCM-NN using daily monitoring data of soybean pod borer collected in Changchun, Jilin Province, from July to September during 2020–2023. Experimental results demonstrate that PCM-NN preserves biological interpretability while exhibiting strong nonlinear representational capacity, offering a feasible pathway for pest modeling and forecasting under multi-factor climatic conditions. This approach provides valuable support for agricultural pest monitoring, prevention, and control strategies.

\end{abstract}


\maketitle
\section{Introduction}
\begin{figure}[h]
\begin{center}
\includegraphics[width=\linewidth]{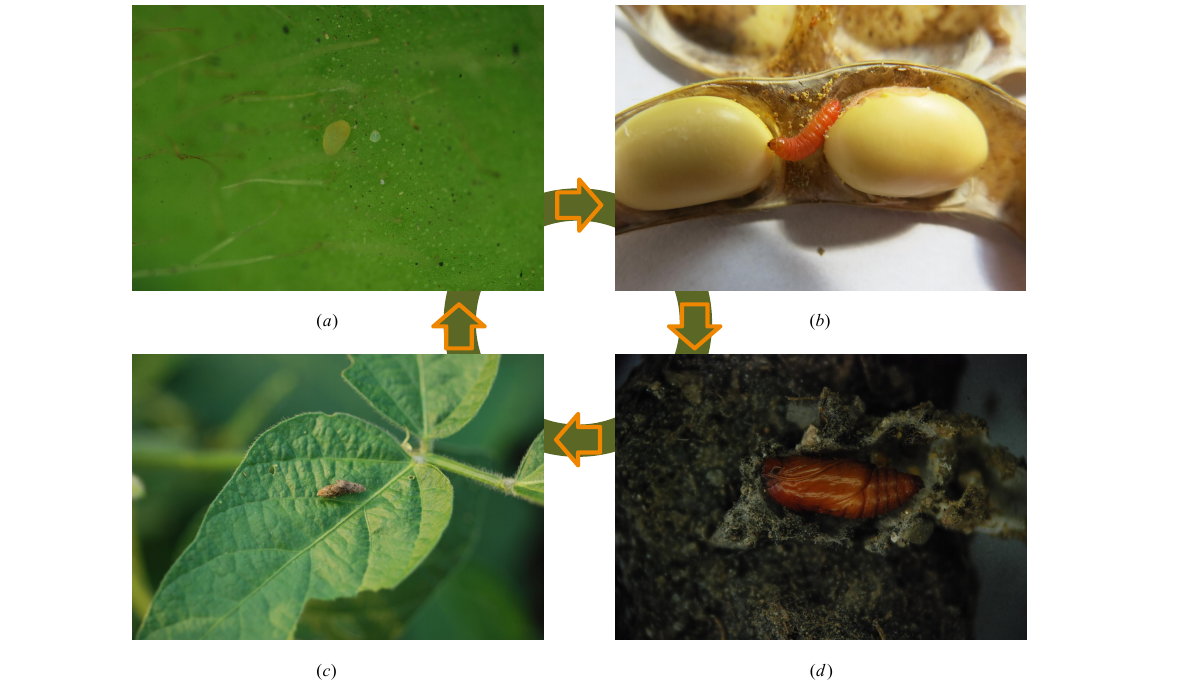}
\end{center}
\caption{ {\bf Different Developmental Stages of the Soybean Pod Borer .} (a) Egg Stage. (b) Larvae Stage. (c) Adult Stage. (d) Pupa Stage. Provided by Prof. Shusen Shi.}
\label{Figure1}
\end{figure}

Soybean is an important grain, oil, and feed crop in China and plays a strategic role in the development of the national economy. However, in recent years, climate change and agricultural ecosystem imbalances have exacerbated pest threats to soybean production, among which the \textit{Leguminivora glycinivorella} (commonly known as the soybean pod borer) is one of the major pests responsible for substantial yield and economic losses. Research indicates that this pest's suitable habitats are widely distributed across most regions of China \cite{gao2018current,shishusen, fei2024integrative}, the Korean Peninsula, the Russian Far East, and Japan \cite{yang2024predicting}. The northeastern region of China, as the country’s primary soybean-producing area, faces a particularly severe threat from soybean pod borer outbreaks, according to the 2023 Technical Plan for the Prevention and Control of Major Soybean Diseases and Pests\footnote{\url{http://www.moa.gov.cn/ztzl/ddymdzfhjs/jszd_29063/202303/t20230301_6421921.htm}}. During the critical soybean growth period (July-September), the larvae bore into pods and feed on seeds, leading not only to direct yield losses but also to a sharp decline in crop market value. For accurate pest identification, Fig. \ref{Figure1} illustrates the typical morphological characteristics of this pest across its life stages. Given the substantial threat posed by soybean pod borers to soybean production, there is an urgent need to conduct modeling and prediction research on the occurrence dynamics of soybean pod borers, so as to achieve early warning and scientific prevention and control of the pest occurrence.

Pest control has long been a pressing global challenge, attracting extensive attention from entomologists and plant protection scientists (\cite{chen2022identifying,jiang2021hyperspectral,kuzmin2020soybean}). Currently, pest prediction research primarily encompasses two categories: empirical models based on statistical correlations and analytical models based on biophysical mechanisms. Statistical approaches establish associations between environmental factors and pest occurrence through regression analysis or machine learning algorithms. For example, \cite{chen2022pest} demonstrated that lychee stink bug (Tessaratoma papillosa) populations thrive under optimal temperature, humidity, and light conditions. \cite{amjad2022effect} investigated the effects of temperature and humidity on cotton pest complexes, quantifying the correlation coefficients between climatic factors and pest abundance. Similarly, \cite{kumar2023seasonal} analyzed the relationship between major peanut pests and microclimate variables. However, while these methods can identify key environmental drivers, they fail to elucidate the physiological mechanisms underlying pest population responses to environmental changes. On the other hand, mechanistic models employ mathematical tools such as differential equations to describe pest population growth and dispersal dynamics  (\cite{sokame2021system,malaguit2023identifying}). Yet, existing models often neglect the dynamic influence of climatic factors or treat them merely as external inputs, rather than integrating temperature, humidity, and other variables as deterministic parameters within the model framework (\cite{skendvzic2021impact,ma2022potential,
ponti2021biological}). A study by \cite{deutsch2018increase} revealed that each 1°C increase in global temperature escalates pest-induced crop losses by 10–25\%, underscoring the urgent need for predictive models that couple climatic variables with biological parameters. Nevertheless, current research still falls short in effectively integrating climatic factors with pest population dynamics, necessitating novel modeling approaches that simultaneously account for climatic influences and biological mechanisms to enhance prediction accuracy and practical applicability.

To overcome the insufficient integration of climatic factors and pest population dynamics in existing modeling approaches, this study proposes a population dynamics model of soybean pod borer coupled with temperature and humidity. Specifically, we formulate a logistic-type ordinary differential equation (ODE) to describe the temporal evolution of pest populations, and embed a microclimate-driven parameter function $\alpha(T, H, t)$ to quantify the combined impact of environmental conditions. To infer the non-explicit function $\alpha$, we integrate the ODE model into a Physics-Informed Neural Network (PINN) framework and develop a novel Pest Correlation Model Neural Network (PCM-NN). By minimizing the composite loss function composed of ODE residuals and observation errors, this model not only fits the observed soybean pod borer population data but also estimates the  evolving form of $\alpha$. Compared to traditional models that require prior assumptions about functional forms, PCM-NN flexibly captures the nonlinear dependencies of $\alpha$ on temperature, humidity, and time, thereby enhancing both the expressiveness and the physical consistency of the model. Furthermore, we utilized monitoring data from 2020 to 2023 to train and validate the model. Through this process, we successfully estimated the important microclimate-driven parameter $\alpha$, which is essential for understanding the mechanisms of pest occurrence and developing effective prevention and control strategies. Additionally, we used the trained model to predict the occurrence of soybean pod borer in 2024, providing forward-looking guidance for agricultural production.


    


This study presents an innovative integration of Physics-Informed Neural Networks (PINNs) with pest population dynamics modeling, resulting in the development of the Pest Correlation Model Neural Network (PCM-NN) for soybean pod borer prediction. By embedding the mechanistic pest dynamics model as a physical constraint within the neural network architecture, PCM-NN achieves accurate fitting of pest population trends while simultaneously quantifying and visualizing the regulatory effects of microclimate factors such as temperature and humidity. Due to its methodological flexibility and generalizability, this approach holds significant promise for broader applications in the prediction of other agricultural pests and provides a feasible technical framework for climate-resilient pest forecasting and precision management.

The remainder of this study is organized as follows. Section 2 introduces the data sources and preprocessing procedures for soybean pod borer monitoring, followed by a detailed description of the climate-coupled logistic model and its implementation within the PINN framework. Section 3 presents the experimental results and validation analyses, assessing the model’s performance and applicability. Section 4 discusses the advantages, limitations, and potential practical implications of the proposed approach. Finally, Section 5 concludes the study and outlines directions for future research.


\begin{figure}[htb]
\begin{center}
\includegraphics[width=\textwidth]{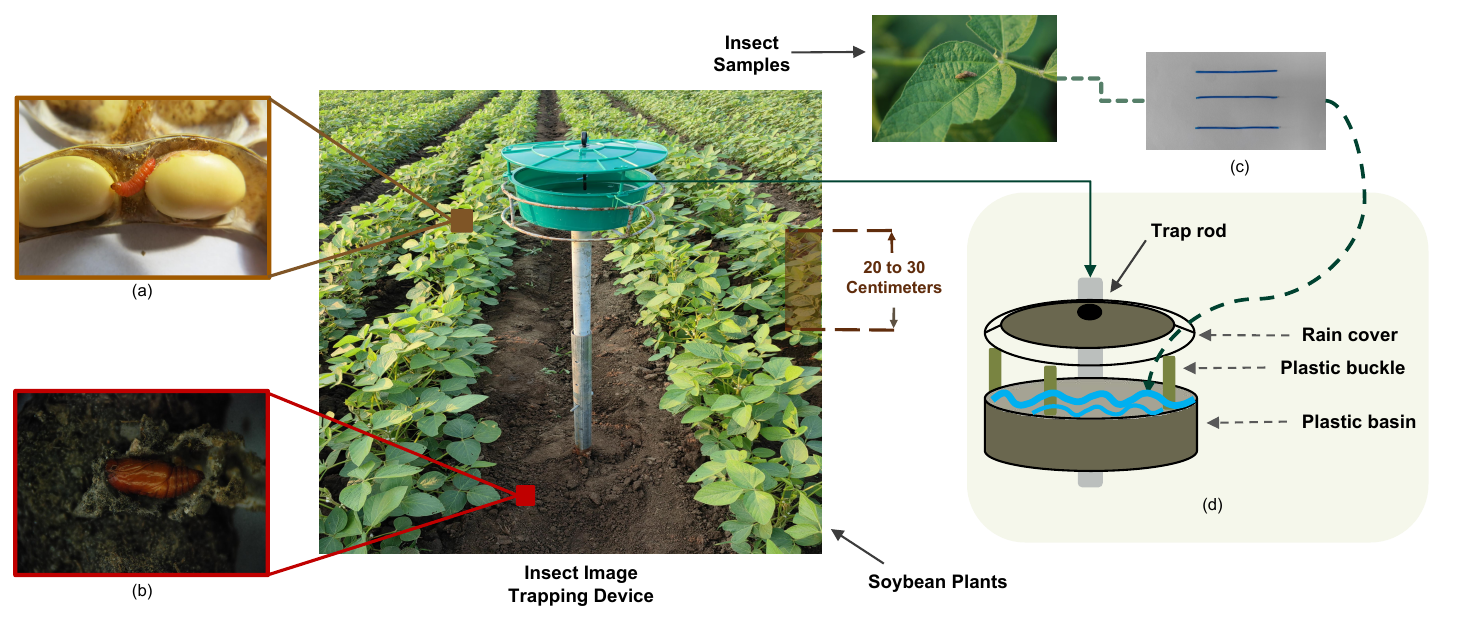}
\end{center}
\caption{{\bf Integrated monitoring system for soybean pod borer (\textit{Leguminivora glycinivorella}) population dynamics.}
(a) Adult male \textit{L. glycinivorella}.
(b) Final instar larva infesting soybean pod.
(c) Synthetic sex pheromone lure.
(d) Pheromone trap (PVC construction) deployed at 20-30 cm above soybean canopy.
Georeferenced deployment shown at Jilin Agricultural University Experimental Station (43.82°N, 125.42°E).}\label{Figure2}
\end{figure}
\section{Data and Methods}
\subsection{Data}
\subsubsection{Data Acquisition}
Deploying lepidopteran sex pheromones constitutes a cornerstone of precision pest surveillance in modern integrated pest management frameworks \cite{wakamura1992development}. For \textit{Leguminivora glycinivorella}, species-specific pheromone components (predominantly (Z)-9-tetradecenyl acetate and (Z)-11-tetradecenyl acetate in optimized ratios) have been chemically characterized and operationally validated for field monitoring \cite{le2006sex}. This semiochemical-based approach exploits male moths' obligate orientation behavior toward synthetic pheromone plumes, enabling real-time tracking of adult eclosion peaks and spatiotemporal population gradients \cite{hu2012synthesis}.

Here, we employ trapping devices to collect data on the male adult soybean pod borer and monitor population dynamics during its peak activity period of emergence. The specific configuration of the trapping device is illustrated in Fig. \ref{Figure2}, and the key settings are listed in the following.
\begin{itemize}
\item {\bf Study site}. Field monitoring was conducted at the Soybean Regional Technology Innovation Center research base (43.82°N, 125.42°E) of Jilin Agricultural University, located in the suburban agroecosystem of Changchun City, Jilin Province, Northeast China. A 1-hectare (10,000 ${\rm m}^2$) continuous soybean cropping field was selected to monitor natural soybean pod borer population dynamics. To eliminate anthropogenic interference, no insecticide applications or pheromone-based mating disruption tactics were implemented during the soybean growing season (May–September).
\item {\bf Trapping devices}. Three traps (T1–T3) were deployed in a triangular grid (inter-trap distance $\geq 50$ m) in the end of June before the anticipated adult emergence period (at July), positioned 20-30 cm above the soybean canopy to intercept male flight paths. In addition, traps were placed $\geq 10$ m from field edges to reduce edge-effect bias. Synthetic sex pheromone lures (Ningbo Newcomb Biotechnology Co., LTd., Ningbo, Zhejiang, China) replaced monthly to maintain emission efficacy. 
\item {\bf A phased observation}. Pre-emergence phase: Tri-daily trap inspections to establish the baseline activity.
Active Monitoring Phase: Upon initial detection of male adults (defined as $\geq 1$ individual/3-trap/day), inspections intensified to daily intervals. 
\item {\bf Meteorological Data}. Daily-averaged microclimatic parameters (mean temperature [°C] and relative humidity [\%]) for Changchun City were obtained from the Changchun National Benchmark Climate Station under the China Meteorological Administration.
\end{itemize}

Integrative analysis of daily trap captures (summed across three spatially distributed pheromone traps) revealed consistent interannual patterns in adult soybean pod borer emergence within the Changchun agroecosystem. The population phenology typically exhibited a unimodal curve or occasional bimodal fluctuations. 
This unimodal/bimodal plasticity likely reflects microclimatic modulation of pupal development synchronicity, with secondary peaks (when observed) correlating with mid-August precipitation events. The standardized trapping system effectively captured these transient dynamics.

\subsubsection{Data Pre-processing}
In this study, daily monitoring records of adult male soybean pod borer populations were systematically collected alongside corresponding mean daily temperature and humidity measurements during July-September annually from 2020 to 2024. Adult emergence is typically initiated in late July (July 25–31), with population peaks occurring between August 8–15, followed by a terminal decline phase in late August that is synchronized with soybean pod hardening (R6 growth stage).

To minimize phenological variability across years and specifically target the period when the interaction between pest dynamics and environmental factors is most relevant for modeling purposes (as the key life-cycle stages of the soybean pod borer and critical soybean growth phases closely align within this timeframe, facilitating the capture of essential ecological relationships), we restricted analyses to a critical phenological phase (July 25–August 23) during which soybean pod borer adults exhibit peak activity in temperate soybean agroecosystems.
The dataset was partitioned chronologically, with observations from 2020 to 2023 allocated to model training and those from 2024 reserved for prospective validation of temporal generalizability.

\begin{figure}[tbp]
\begin{center}
\includegraphics[width=\linewidth]{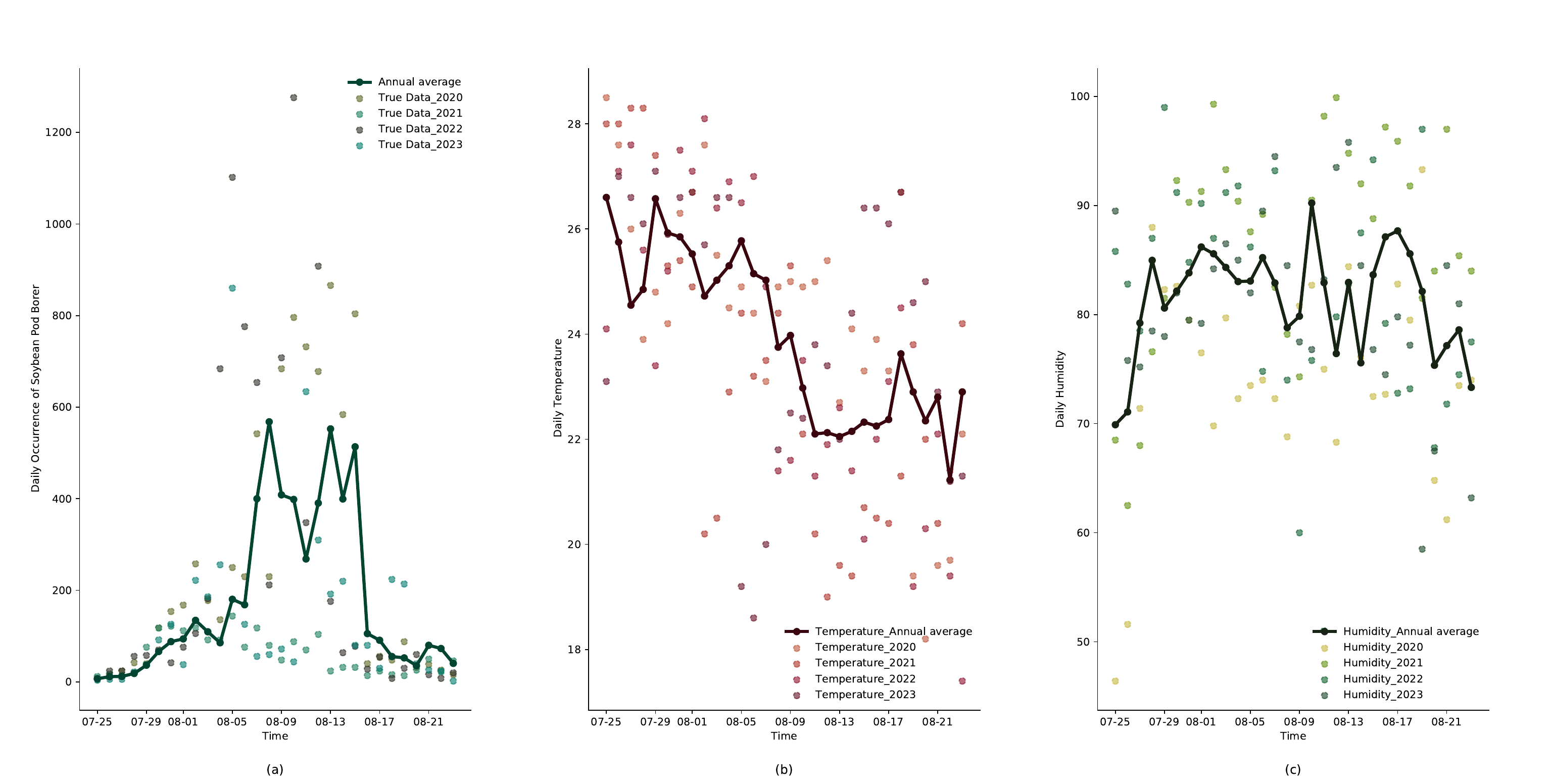}
\end{center}
\caption{{\bf Spatiotemporal dynamics of soybean pod borer (\textit{Leguminivora glycinivorella}) abundance and associated microclimatic drivers during peak activity seasons (2020--2023).}
(a) Daily male moth captures from July 25 to August 23, with individual annual trajectories (colored points) against the interannual average (colored line).
(b) Corresponding daily mean temperature (°C) profiles, illustrating thermal variability across years (individual observations) and climatic baseline (interannual average).
(c) Daily relative humidity (\%) fluctuations (individual observations and interannual average) during the same period.
Interannual averages derived from arithmetic means of synchronized daily records across 2020--2023. Microclimatic data sourced from China Meteorological Administration.}
\label{Figure3}
\end{figure}
Initial comparative analysis of daily observations across consecutive years (2020–2023) revealed pronounced instability and irregularity in 
soybean pod borer population counts and their associated microclimatic variables (temperature and humidity). These fluctuations, indicative of high stochastic noise inherent in short-term ecological datasets, obscured discernible phenological or climatic trends. To mitigate this limitation, we adopted interannual averaging, a climatologically inspired preprocessing method, to disentangle persistent seasonal signals from transient perturbations, as shown in Fig. \ref{Figure3}. It should be noted in Fig. \ref{Figure3} that, given the roughly equal proportion of male and female adult soybean pod moths in the field population, the data shown in Fig. \ref{Figure3} effectively represents the number of adult populations by multiplying the collected number of males by two.

Interannual averaging leverages temporal compositing to suppress year-specific anomalies (e.g., extreme weather events, localized management practices) while amplifying recurrent patterns driven by climatic seasonality. This approach aligns with methodologies used in climate science to define baseline bioclimatic norms \cite{RN4} and has been adapted here for agroecological systems to resolve insect-environment coupling. The procedure entailed two sequential steps:
\begin{enumerate}
    \item Temporal Alignment: For each calendar date between July 25 and August 23 (Day of Year [DOY]), we aggregated four years (2020–2023) of synchronized daily records, including: adult male soybean pod borer counts (individuals/3-trap/day), mean daily temperature (°C), mean relative humidity (\%).
    \item Composite Dataset Generation:
The interannual composite was derived by computing the arithmetic mean for each variable at every DOY:
$$\bar{x}_{DOY}=\frac{1}{N}\sum_{i=2020}^{2023}x^{i}_{DOY},$$
where $\bar{x}_{DOY}$ represents the averaged value for a given day-of-year, and $N=4$ (years). This produced a 30-day synthetic dataset reflecting the "climatological normal" of pest incidence and microclimate during peak activity seasons.
\end{enumerate}
To contextualize 2024 test data against historical patterns, an expanded composite (2020–2024) was generated using identical procedures, though model training strictly utilized the 2020–2023 subset to prevent temporal leakage.

\begin{table}[h]
\centering
\caption{\bf Summary of Key Studies on Crop Disease and Pest Analysis Using Logistics}
\label{Table1}
\renewcommand{\arraystretch}{2}
\begin{tabular}{p{1.5cm}p{4cm}p{4.5cm}p{4.5cm}}
\hline
Reference & Methodology & Key Findings & Limitations \\
\hline
\cite{tab1-1} & Math models for analyzing optimal harvesting policies & Unique periodic solutions and optimal harvesting efforts identified & Focuses mainly on theoretical models, practical application details lacking \\
\cite{tab1-2} & Math analysis of cumulative lethal effect of pesticide spraying & Conditions for pest extinction and stability of periodic solutions discovered & Theoretical research; experimental validation and practical application details limited \\
\cite{tab1-3} & Logistic models and bioassays for insecticide effectiveness & Factors influencing cockroach mortality identified; application frequencies recommended & Limitations not specified in excerpt \\
\cite{tab1-4} & Controlled experiments and statistical models for fungi effectiveness & Soil moisture and temperature impact on fungi virulence and insect susceptibility found & Specific fungal isolates and insect species tested; controlled conditions may not reflect natural variations \\
\cite{tab1-5} & Math models and simulations for integrated pest control & Complex dynamical behaviors of pest-parasitoid populations; periodic control strategies suggested & Simplifying assumptions may limit applicability to real-world scenarios \\
\cite{tab1-6} & Logistic reaction-diffusion harvesting model with infinite delay & Evolution rate and small diffusion rate benefit species survival & Simplifies complex ecological phenomena; further empirical validation needed \\
\hline
\end{tabular}
\end{table}
\subsection{Methods}
\subsubsection{Model Validation}

The logistic equation, a classic mathematical model for describing population growth patterns, was first proposed by Verhulst in the mid-19th century to explore the relationship between population growth and resource constraints. Let $x(t)$ represent the population size at time $t$, the classical Logistic equation is
\begin{align}\label{equation1}
    \frac{\mathrm{d}x}{\mathrm{d}t}= rx\left(1-\frac{x}{K}\right),
\end{align}
where $r$ represents the intrinsic growth rate, and $K$ represents the carrying capacity of the environment. Here, the first term $rx$ is used to model unimpeded growth rate in the beginning period; and the second term $-\frac{r}{K}x^2$ measures the interference between population by competing for some critical resource, such as food or living space.
Due to its simplicity and generality, the logistic equation has become a fundamental tool in ecology and complex systems science.
Many researchers integrate this differential equation into pest population dynamics to predict development trends and optimize control strategies.
To gain a deeper understanding of the application of this model in crop pest research, we conducted an extensive literature review and summarized the relevant research findings in Table \ref{Table1}. 

\begin{figure}[htb]
\begin{center}
\includegraphics[width=\linewidth]{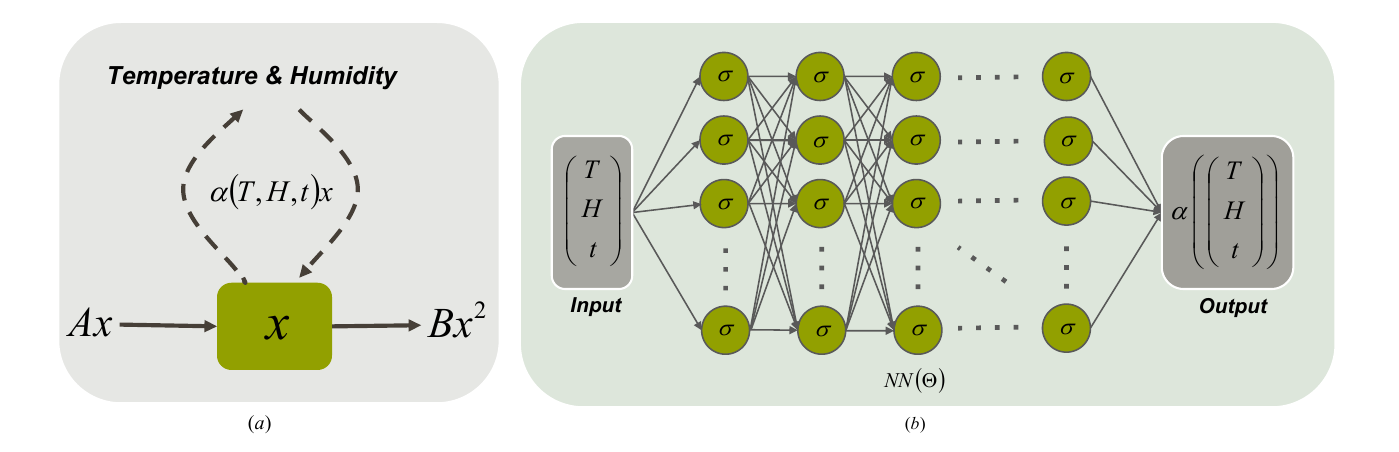}
\end{center}
\caption{(a) Mechanistic-climate feedback structure of Model \eqref{equation2}. Solid black arrows denote density-dependent growth/decline processes, dashed gray arrows represent microclimatic couplings. (b) A deep neural network architecture for climatic modulation function $\alpha^{NN}$.}\label{Figure4}
\end{figure}

The classical logistic model provides a phenomenological framework for density-dependent population growth, yet its assumption of constant intrinsic growth rate $r$ and static carrying capacity $K$ limits applicability to poikilothermic organisms like Leguminivora glycinivorella, whose vital rates are thermally modulated. To mechanistically couple microclimatic factors with pest phenology, we propose a modified logistic equation:
\begin{align}\label{equation2}
\frac{\mathrm{d}x}{\mathrm{d}t} = \underbrace{\left[A + \alpha(T, H, t)\right]x}_{\text{Bidirectional climatic modulation}} - \underbrace{Bx^2.}_{\text{Density limitation}}
\end{align}
Here, $x(t)\in \mathbb{R}^+$ represents the adult population of soybean pod borer at time $t$; $A\in \mathbb{R}^+$ is the field-realistic baseline growth rate, reflects growth potential under average field conditions, which is inherently lower than laboratory-observed maxima due to persistent ecological pressures (e.g., predation, and suboptimal microhabitats); $B=\frac{A}{K}$ represents the density-dependent inhibition coefficient, where $K$ is the environmental carrying capacity. In addition, $\alpha: {\mathbb{R}^+}^3 \rightarrow \mathbb{R}$ is the climatic modulation function defined as:
\begin{equation}\label{equation3}
\alpha(T, H, t) = \alpha\left(\left(T(t)-T^*\right)^2, \left(H(t)-H^*\right)^2, t\right)
\end{equation}
with $T^*$ and $H^*$ denoting thermal and hygric optima for the survival of soybean pod borer. 
A possible property of bidirectional climatic modulation function $\alpha$ is
\begin{equation*}
\left\{
\begin{aligned}
&\alpha(T, H, t)\geq 0, ~~~\text{(growth facilitation)} \qquad &\text{if} \left(T(t)-T^*\right)^2\leq M_1~\text{and}~\left(H(t)-H^*\right)^2\leq M_2,\\
&\alpha(T, H, t)< 0, ~\text{(growth suppression)} \quad 
&\text{if} \left(T(t)-T^*\right)^2> M_1~\text{or}~\left(H(t)-H^*\right)^2> M_2.
\end{aligned}
\right.
\end{equation*}
quantifying how Leguminivora glycinivorella populations experience growth suppression under thermal/hygric stress (deviations exceeding thresholds $M_1, M_2$) versus growth facilitation when microclimatic conditions remain within optimal ranges.
However, the functional form of the climatic modulation term inherently encapsulates complex, non-clear interactions among temperature, humidity, and pest physiology. 
Traditional explicit parameterizations (e.g., quadratic penalties) risk oversimplifying these higher-order couplings, our enhanced framework employs a fully connected neural network (NN) to represent $\alpha(T, H, t)$, enabling data-driven discovery of nonlinear interactions between microclimatic factors and pest population dynamics. 
The transfer diagram of model \eqref{equation2} is shown in Fig. \ref{Figure4}(a), and the NN diagrammatic representation of \eqref{equation3} is given in Fig. \ref{Figure4}(b).
\subsubsection{PCM-NN}
This section introduces a new deep learning method - Pest Correlation Model Neural Network (PCM-NN). This method not only inherits the essence of Physical Information Neural Network (PINN), which improves the accuracy of model prediction by incorporating physical laws, but also introduces additional environmental parameters (temperature and humidity) to meet the practical needs of pest management, thereby achieving more accurate description and prediction of pest dynamics.

\begin{figure}[t]
\begin{center}
\includegraphics[width=\linewidth]{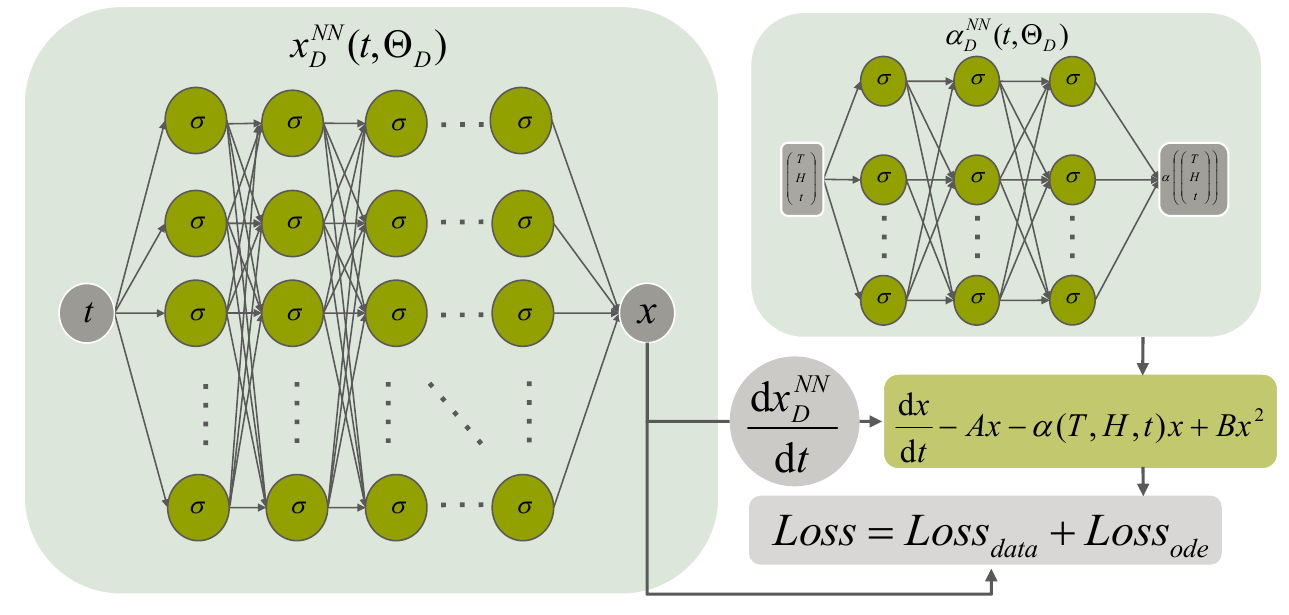}
\end{center}
\caption{{\bf Schematic diagram of PCM-NN algorithm.}
$x_D^{NN}(t,\Theta_{D})$ serves as the state variable for fitting Model \eqref{equation3}, as depicted in the top-left region of the figure. 
Concurrently, parameter $\alpha_D^{NN}(t,\Theta_{D})$ is utilized for inferring the parameter $\alpha$, illustrated in the top-right region. 
The term $\frac{\mathrm{d}x_{D}^{NN}}{\mathrm{d}t}$ denotes the self-differential operator, central to the model's dynamics. 
The overall Loss function, $Loss$, comprises two components: $Loss_{data}$ and $Loss_{ode}$. 
By minimizing $Loss$, the algorithm achieves simultaneous data fitting and inference of parameter.
}\label{Figure5}
\end{figure}

Fig. \ref{Figure4} (b) shows a schematic diagram of the deep neural network structure of PCM-NN. Unlike traditional PINNs that only receive time $t$ as input, PCM-NN's input layer is expanded to include three dimensions: time $t$, ambient temperature $T$, and ambient humidity $H$. 
This design aims to capture the impact of environmental factors on the occurrence and development of pest infestations, thereby improving the predictive accuracy of the model.

Fig. \ref{Figure5} provides a detailed description of the algorithm flow of PCM-NN. The core lies in embedding the pest model \eqref{equation2} as physical information into the training process of the neural network. This pest model is based on ecological principles and describes the variation of pest numbers over time, temperature, and humidity. In PCM-NN, this physical information is represented by defining a specific Loss function, which includes two parts:
\begin{itemize}

\item {\bf Data Loss}: measures the difference between network predictions and observed pest data. This Loss ensures that the model can learn the dynamic behavior of pests from historical data. The specific form of Loss is:
\begin{equation}
\begin{aligned}
Loss_{data} & =\frac{1}{N_d} \sum_{i=1}^{N_d}\left|x^{N N}\left(t_i\right)-x\left(t_i\right)\right|^2, 
\end{aligned}
\end{equation}
where $N_d $ represents the quantity of training data, and $x^{NN}(t_{i})$ represents the daily increase in the number of soybean pod borer.

\item {\bf ODE (Ordinary Differential Equation) Loss}: Based on the physical constraints defined by the pest model (2), calculate the deviation of network predictions from physical laws. This Loss term ensures that the model predictions not only conform to observed data, but also strictly follow the physical mechanisms of pest occurrence. 
The form of $Loss_{ode}$ is as follows:
\begin{equation}
Loss_{ode}=\frac{1}{N_e}\sum_{i=1}^{N_e}\left|L(t_{i})\right|^2, 
\end{equation}
where $N_e$  denotes the quantity of residual points, chosen at random from the entire computational domain and 
\begin{equation}
\begin{aligned}
L\left(t_i\right)=\frac{\mathrm{d}x^{NN}\left(t_i\right)}{\mathrm{d} t} -\left(Ax^{NN}(t_{i})+\alpha^{NN} \left(T\left(t_i\right), H\left(t_i\right), t_i\right)x^{NN}(t_{i})-B\left( x^{{NN}}\left(t_i\right)\right)^2\right).
\end{aligned}
\end{equation}
\end{itemize}

The total Loss function of PCM-NN is the weighted sum of $Loss_{data}$ and $Loss_{ode}$, and the choice of weights depends on the specific requirements and data characteristics of the problem. By minimizing this total Loss, PCM-NN can simultaneously utilize observational data and physical laws to achieve efficient and accurate pest prediction. 

\subsection{Experiment Settings}
\begin{figure}[htb]
\begin{center}
\includegraphics[width=\linewidth]{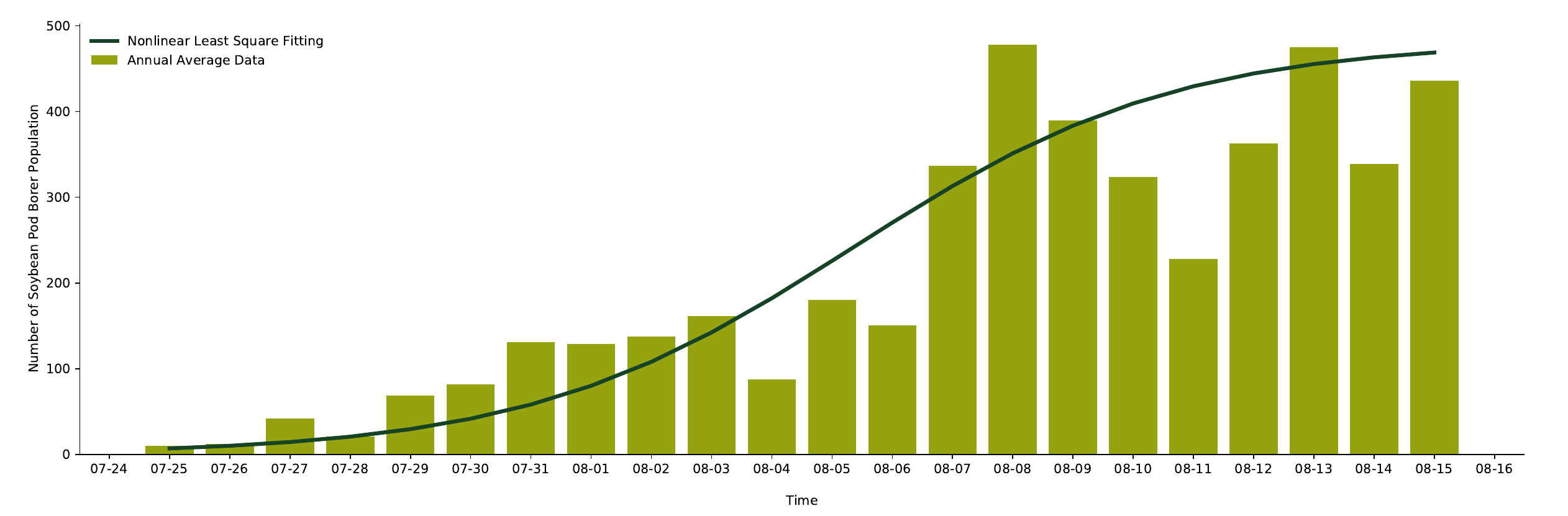}
\end{center}
\caption{{\bf Nonlinear least squares fitting results.} 
The bar chart in the figure represents the real data of soybean pod borer. 
The line (represented in green) is a curve obtained by fitting these real data points using nonlinear least squares method.
 }
\label{Figure6}
\end{figure}

Firstly, we determined the values of parameters $A$ and $B$. 
Utilizing nonlinear least squares, we employed the annual seasonal average data of soybean pod borer from 2020 to 2023, along with the corresponding temperature and humidity data, and selected the first 1-22 days for nonlinear least squares fitting. 
The reason for making this choice is that, based on data observation, the initial 22 day data period is relatively close to the trend of the logistic growth curve.
Consequently, we obtained the parameter values of $A=0.372$ and $B=0.0008$. 
Fig. \ref{Figure6} presents the fitting results of the nonlinear least squares analysis.
Our experimental setup aimed to ensure that the neural network model achieves optimal performance in predicting soybean pod borer occurrences, while taking into account environmental factors and specific characteristics of the pest.

As illustrated in Table \ref{Table2}, the depth and width of the neural network employed in our study, along with the number of iterations for the algorithm, were configured to optimize model performance. 
These configurations encompass the number of hidden layers, the quantity of neurons per layer, and the selection of activation functions, all of which were based on experimental validation and performance evaluation results.
Crucially, the table also lists the optimal temperature and humidity for the soybean pod borer, which are 21°C and 84\%, respectively. These values are derived from the fundamental survival characteristics of the pest itself. 
Additionally, we constrained the parameter $\alpha$ within the range of $[-1.372,0.628]$, ensuring that the range of $A-\alpha$ remains within $[-1,1]$.

\begin{table}[thb]
 \centering
\caption{\bf Deep Neural Network Hyperparameter Values and Parameter Settings}
\label{Table2}
\renewcommand{\arraystretch}{1.5} 
\begin{tabular}{p{4cm}p{9cm}p{2cm}}
\hline
Symbol & Parameter Interpretation & Values 
\\ \hline
$T^*$              & Optimum temperature                                  & $21^\circ C$ \\
$H^*$              & Optimum humidity                                     & 84\% \\
$A$                & Field-realistic baseline growth rate                         & 0.372\\
$B$                & Field-realistic baseline growth rate/carrying capacity & 0.0008\\
$\alpha$           & Climatic modulation function                                  & [-1.372, 0.628]\\
$U^{NN}(t)$        & Number of hidden layers in data                      & 5*32\\
$\alpha^{NN}(t)$   & Number of hidden layers for parameters               & 3*64\\
Iterations         & Number of iterations                                 & $1 \times 10^4$
\\ \hline
\end{tabular}
\end{table}
\section{Result}
\subsection{Fitting Results and Inferred Parameter}
\begin{figure}[htb]
\begin{center}
\includegraphics[width=\linewidth]{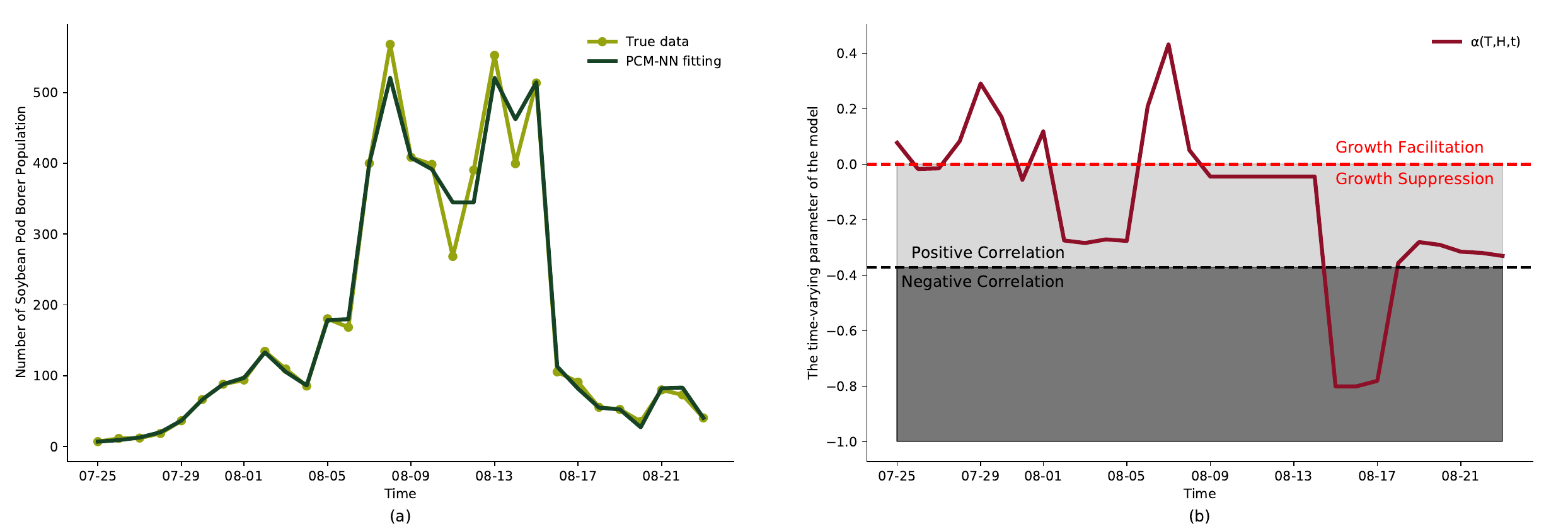}
\end{center}
\caption{{\bf PCM-NN fitting inference results.} 
(a) The fitting results of the PCM-NN algorithm are presented. The light green line represents the actual data, while the dark green line depicts the data fitted by the PCM-NN algorithm.
(b) The estimated value of the parameter $\alpha$ inferred from the fitting results of PCM-NN is shown as a parameter curve over time, which is represented by a red line. The red dashed line represents the growth threshold line, while the black dashed line represents the un-correlated baseline.
}
\label{Figure7}
\end{figure}

This subsection presents the fitting results of our proposed PCM-NN algorithm. By leveraging the PCM-NN, we aim to capture the dynamic nature of soybean pod borer outbreaks more accurately than traditional methods.
To illustrate the fitting effects, Fig. \ref{Figure7} presents two graphs. 
Fig. \ref{Figure7} (a) shows the comparison between the actual outbreak of soybean pod borer and the fitted value of PCM-NN algorithm 
This graph clearly demonstrates the model's ability to capture the trends and fluctuations in the outbreak data, highlighting its robustness and reliability.

More importantly, Fig. \ref{Figure7} (b) showcases the inferred parameter $\alpha$ over time. Unlike traditional estimation methods that provide a single, fixed value for alpha, our PCM-NN algorithm captures its dynamic nature. 
This parameter $\alpha$ provides a more realistic representation of the underlying processes governing soybean pod borer outbreaks, reflecting changes in environmental conditions, pest management practices, and other relevant factors. 

More importantly, Fig. \ref{Figure7} (b) showcases the parameter $\alpha$ inferred over time by our PCM-NN algorithm. 
Unlike traditional estimation methods that provide a single, fixed value for $\alpha$, our approach captures its dynamic nature. 
This parameter $\alpha$ offers a more realistic representation of the underlying processes governing soybean pod borer outbreaks, reflecting variations in environmental conditions, pest management practices, and other pertinent factors.

Upon observation, it is evident that the trend of $\alpha$ relatively aligns with that of the actual data. 
Specifically, $\alpha$ tends to increase when the population quantity (or other relevant metrics) exhibits an upward trend. 
Conversely, during periods of data decline, $\alpha$ also exhibits a downward trend. 
While this correlation is not perfectly synchronous due to potential time lags inherent in the system, it is noteworthy that when $\alpha$ is greater than zero, it indicates growth promotion, whereas when $\alpha$ is less than zero, it indicates growth suppression. Furthermore, when $\alpha$ exceeds -0.372, it exhibits a positive correlation with the data, whereas when $\alpha$ falls below -0.372, it shows a negative correlation. 
This nuanced understanding of the dynamic behavior of $\alpha$ provides deeper insights into the mechanisms driving soybean pod borer population dynamics.

\subsection{Verification by ODE Solution and Predictive Performance of PCM-NN}
To further validate the training effectiveness of our PCM-NN algorithm, we conducted two key tests: ODE back-solving verification and predictive performance assessment.
\begin{figure}[thb]
\begin{center}
\includegraphics[width=\linewidth]{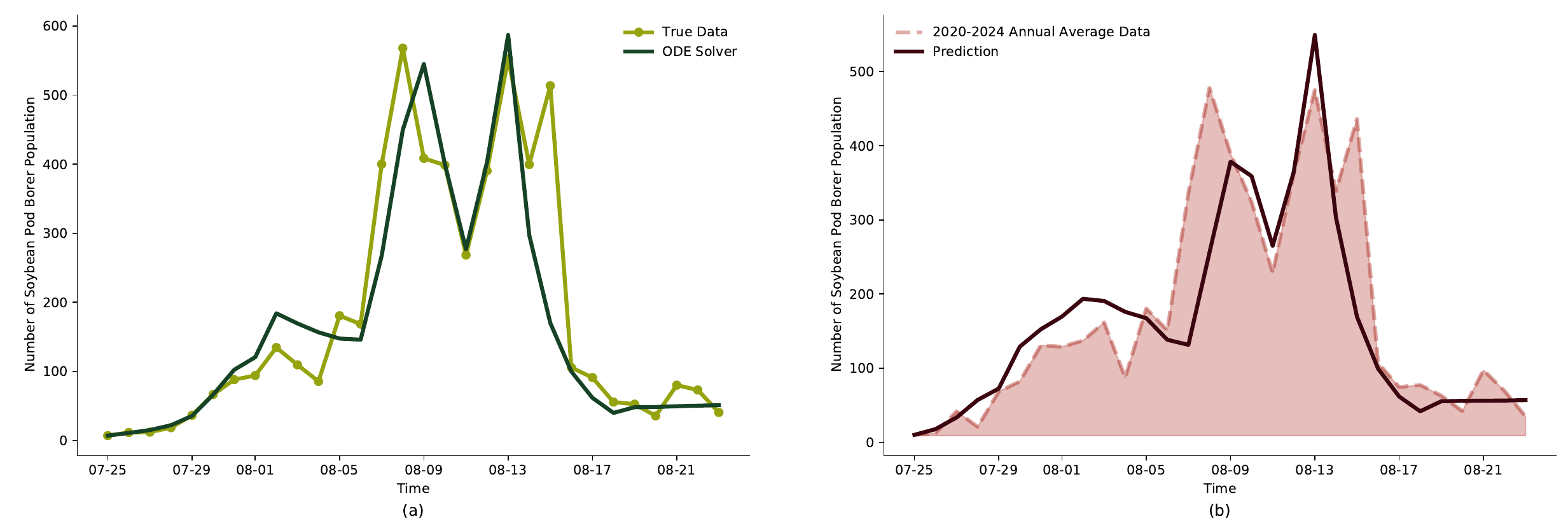}
\end{center}
\caption{{\bf The ODE inverse solution and prediction results of PCM-NN.} (a) Comparison between real data and ODE inverse solution. The light green line represents the actual data, while the dark green line represents the ODE inverse solution. (b) Comparison between forecast curve and interannual average data including 2024. The red dashed line represents the actual data, while the black solid line represents the predicted curve.}
\label{Figure8}
\end{figure}
\subsubsection{Verification}
We verified the trained PCM-NN by substituting the inferred parameter $\alpha$ back into the modified logistic equation \eqref{equation2} and solving it using a ODE solver. 
Fig. \ref{Figure8} (a) presents the results of this ODE back-solving. 
The simulated curve, obtained from the pest model with the inferred climatic modulation function 
$\alpha$, aligns well with the observed data. 
This alignment suggests that the training outcomes of PCM-NN are consistent with the underlying transmission mechanisms of the pest model. 
It confirms that PCM-NN has accurately captured the evolution of pest populations, as reflected in the changing nature of $\alpha$.
\subsubsection{Predictive Performance Assessment}
We evaluated the predictive capability of PCM-NN in forecasting the daily occurrence of soybean pod borer during July and September 2024. Utilizing the observed data from the first day of 2024 annual average data set, we input this initial value into the model \eqref{equation2}, incorporating the estimated $\alpha$ by PCM-NN, to predict the subsequent occurrence trend. 
The prediction results are depicted in Fig. \ref{Figure8} (b). 
The predictive curve in Fig. \ref{Figure8} (b) exhibits a similar trend to both the interannual average line spanning from 2020 to 2024. It is worth noting that the predicted values also learned the trend of double peaks and are consistent with the real data.
This observation underscores PCM-NN's ability to generate reliable predictions based on the inferred $\alpha$ and the pest model. 
It demonstrates that PCM-NN not only fits historical data accurately but also possesses the predictive power to anticipate future pest infestation trends. 
The ODE back-solving verification confirms the consistency of PCM-NN's training results with the pest model's mechanisms. 
Additionally, the predictive performance assessment highlights PCM-NN's capacity to forecast future soybean pod borer infestation trends with reasonable accuracy. 
These findings collectively demonstrate the robustness and practicality of PCM-NN in addressing real-world pest management challenges.
\subsection{Fitting Accuracy on Three Evaluation Metrics}
\begin{figure}[t]
\begin{center}
\includegraphics[width=\linewidth]{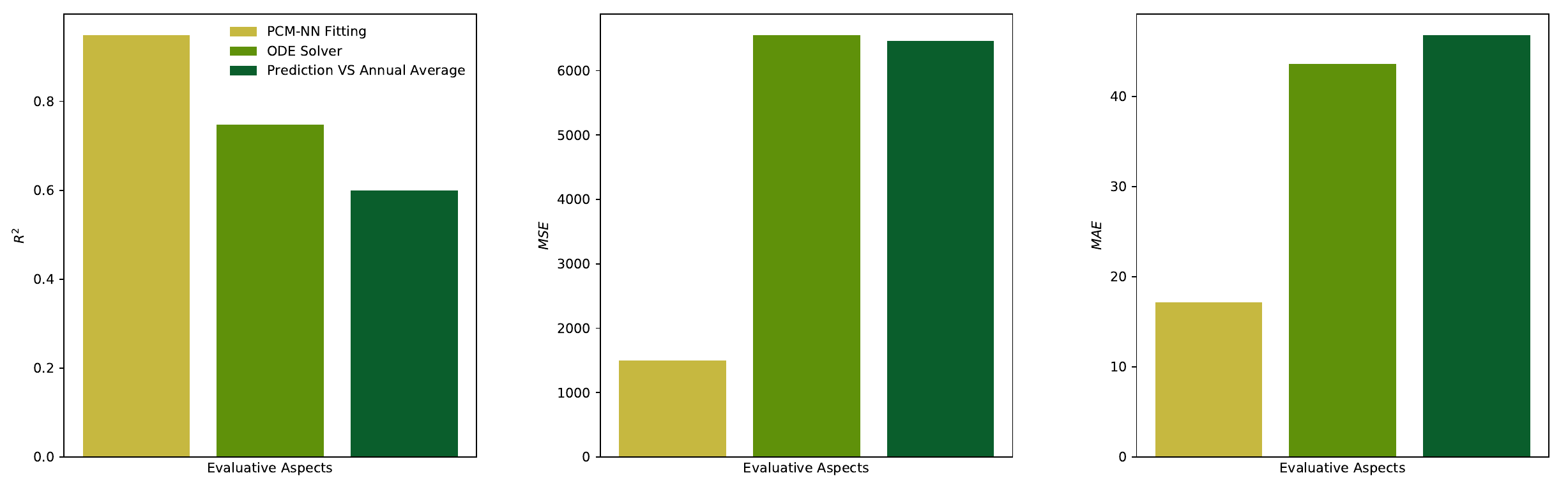}
\end{center}
\caption{{\bf Comparison of \textit{MSE}, \textit{MAE}, and $R^2$ Values for Different Evaluative Aspects.}
}
\label{Figure9}
\end{figure}

The PCM-NN algorithm demonstrates excellent fitting performance, as evidenced by the three key metrics: Mean Squared Error (\textit{MSE}), Mean Absolute Error (\textit{MAE}) and the coefficient of determination ($R^2$). 
Specifically, the \textit{MSE} value indicates the average of the squares of the errors—that is, the average squared difference between the estimated values and the actual data points. A lower \textit{MSE} value suggests a better fit. 
\textit{MAE} measures the average absolute difference between predicted values and actual observed values.
The calculation method is to take the absolute difference between the predicted value and the true value of each sample, and then calculate the average of all samples.
Similarly, the $R^2$ value measures the proportion of the variance in the dependent variable that is predictable from the independent variables. An $R^2$ value closer to 1 indicates a stronger relationship between the predicted and actual values. The formulas for \textit{MSE}, \textit{MAE} and $R^2$ are as follows: 
$$R^2 = 1 - \frac{\sum_{i=1}^{n} (y_i - \hat{y}_i)^2}{\sum_{i=1}^{n} (y_i - \bar{y})^2}, \quad \textit{MSE}=\frac{1}{n}\sum^n_{i=1}(y_{i}-\hat{y}_{i})^2,\quad
\textit{MAE}=\frac{1}{n}\sum^n_{i=1}(y_{i}-\hat{y}_{i}).$$ 
Where $n$ is the total number of samples, and $y_i$ represents the true value of the $i$-th sample, which is the actual observed data. $\hat{y}_i$ indicate the predicted value of the $i$-th sample. 
As shown in Fig. \ref{Figure9}, PCM-NN demonstrates excellent fitting performance across three metrics, indicating that the model fits the data very well. Furthermore, when the parameter $\alpha$ inferred using PCM-NN is substituted into the pest population model \eqref{equation2} and solved using an ODE solver, the resulting data still shows a high level of consistency with the real data, confirming the accuracy of the inferred parameter $\alpha$. Additionally, when compared with the interannual average data from 2020 to 2024, the model’s predictions show a high degree of agreement with the true values, further validating the reliability of our proposed method.
\section{Discussion}

In this study, we developed the pest correlation model neural network (PCM-NN) algorithm, an improvement upon the physics-informed neural network (PINN) framework, specifically designed for fitting and predicting the soybean pod borer population dynamics. 
The PCM-NN framework mitigates the inherent limitations of traditional forecasting methods by integrating a mechanistic constraint differential equation model derived from pest physiology and environmental interactions. 
It reduces reliance on purely data-driven assumptions while maintaining adaptability to nonlinear ecological dynamics.
The well-known classical Logistic model could effectively capture three-phase population dynamics: lag (slow growth) phase,  exponential growth phase, and the stationary phase. By incorporating a couple of microclimatic factors (temperature $T$ and humidity $H$) through the bidirectional climatic modulation parameter function $\alpha(T, H, t)$, we resolve the classical model's oversimplification of environmental modulation. Diverging from traditional explicit parameterizations of $\alpha$, our enhanced framework employs a fully connected neural network (NN) to represent it. We employ an enhanced PINN architecture to infer time-dependent parameter function $\alpha$, enabling data-driven discovery of nonlinear interactions between microclimatic factors and pest population dynamics while preserving ecological interpretability via embedded biophysical constraints.   

A standardized field monitoring was implemented to synchronize pest population dynamics with microclimatic factors across five consecutive growing seasons from 2020 to 2024 at the Soybean Regional Technology Innovation Center research base of Jilin Agricultural University, located in the suburban agroecosystem of Changchun City, Jilin Province, Northeast China. Here, we employ sex pheromone-baited traps to collect data on the male adult soybean pod borer and monitor population dynamics during the soybean growing season (May-September); microclimate data comprising daily-averaged temperature ($^\circ$C) and relative humidity (\%) records were obtained from Changchun National Benchmark Climate Station under the China Meteorological Administration. We implemented an interannual compositing processing to mitigate high-frequency observational noise inherent in time-series data, preserving critical phenological signals for mechanistic modeling. Through the fitting results of the PCM-NN algorithm, we systematically evaluated its predictive performance and predicted future occurrence trends of soybean pod borers.
These research findings deepen the theoretical understanding of soybean pod borer ecology and provide practical tools for precision agriculture. 

Although the PCM-NN algorithm demonstrates strong performance in both fitting and prediction tasks, certain limitations remain. For instance, its sensitivity to data noise may introduce bias into prediction outcomes. Future work should focus on optimizing the algorithm to enhance its robustness against noisy or incomplete data. Moreover, conducting sensitivity analyses of model parameters is essential to ensure the model's stability and reliability across diverse environmental conditions. Techniques such as regularization, data augmentation, and advanced noise-handling mechanisms can be employed to further improve the model’s generalization capability and predictive accuracy.
Currently, the model primarily considers temperature and humidity as environmental factors influencing the population dynamics of soybean pod borers. However, pest outbreaks are typically driven by the interplay of multiple factors. Future studies should incorporate additional variables, such as soil type, crop cultivar, and interspecies biological interactions, to construct a more comprehensive, multi-factor coupled model. Such an approach would offer deeper insights into the complex ecological mechanisms of pest emergence and provide a more robust basis for precision pest control strategies.
While the PCM-NN algorithm performs well in long-term forecasting, its predictive accuracy may degrade over extended periods. Therefore, its integration into adaptive pest management systems is critical. By continuously monitoring pest population trends and environmental changes and regularly updating model parameters and predictions, real-time forecasting and dynamic responses to pest outbreaks can be achieved. This adaptive strategy would significantly enhance the efficiency and efficacy of pest control measures.
Future research should conduct cross-regional validations across diverse geographical areas and climatic conditions to evaluate the PCM-NN algorithm's generalizability and scalability. Comparing model predictions with actual pest occurrences in different locations will help assess its adaptability and accuracy. Furthermore, integrating the PCM-NN algorithm with complementary pest control technologies, such as biological and chemical control, could lead to the development of a more holistic and effective pest management system, thereby offering broader technical support for sustainable agricultural production.

\section{Conclusions}

This study establishes a hybrid neural-mechanistic framework (PCM-NN) that resolves temperature-humidity-pest couplings in soybean pod borer population dynamics. Validated against 2020-2023 field data from Changchun, China, PCM-NN achieves high prediction accuracy for 2024 outbreaks while maintaining ecological interpretability through embedded biological constraints. The model's capacity to balance data-driven learning with pest physiology principles offers actionable intelligence for precision pesticide deployment. Future enhancements will optimize data collection methods, improve data quality, and explore additional influencing factors.

\printcredits
\section*{Funding}
Research of Suli Liu is supported by the National Natural Science Foundation of China (12301627), the Science and Technology Research Projects of the Education Office of Jilin Province, China (JJKH20250046KJ), the Technology Development Program of Jilin Province, China (20210508024RQ). Research of Yu Gao is supported by the Earmarked Fund for China Agriculture Research System of MOF and MARA (Grant No. CARS04), the National Key Research and Development Program of China (2023YFD1401000). 
\section*{Declaration of competing interest}
The authors declare that they have no known competing financial interests or personal relationships that could have appeared to influence the work reported in this paper.
\section*{Data availability}
Data used in the study, such as weather data and pest data, are available upon request to the corresponding author. The source code in Python is available at the following repository: \href{https://github.com/jluorganization/PCM-NN}{https://github.com/jluorganization/PCM-NN}.

\bibliographystyle{cas-model2-names}
\bibliography{cas-refs}

\begin{thebibliography}{26}
\expandafter\ifx\csname natexlab\endcsname\relax\def\natexlab#1{#1}\fi
\providecommand{\url}[1]{\texttt{#1}}
\providecommand{\href}[2]{#2}
\providecommand{\path}[1]{#1}
\providecommand{\DOIprefix}{doi:}
\providecommand{\ArXivprefix}{arXiv:}
\providecommand{\URLprefix}{URL: }
\providecommand{\Pubmedprefix}{pmid:}
\providecommand{\doi}[1]{\href{http://dx.doi.org/#1}{\path{#1}}}
\providecommand{\Pubmed}[1]{\href{pmid:#1}{\path{#1}}}
\providecommand{\bibinfo}[2]{#2}
\ifx\xfnm\relax \def\xfnm[#1]{\unskip,\space#1}\fi
\bibitem[{Amjad~Bashir et~al.(2022)Amjad~Bashir, Batool, Khan, Shahid~Nisar, Farooq, Hashem, Alamri, A.~El-Zohri, Alajmi, Tahir et~al.}]{amjad2022effect}
\bibinfo{author}{Amjad~Bashir, M.}, \bibinfo{author}{Batool, M.}, \bibinfo{author}{Khan, H.}, \bibinfo{author}{Shahid~Nisar, M.}, \bibinfo{author}{Farooq, H.}, \bibinfo{author}{Hashem, M.}, \bibinfo{author}{Alamri, S.}, \bibinfo{author}{A.~El-Zohri, M.}, \bibinfo{author}{Alajmi, R.A.}, \bibinfo{author}{Tahir, M.}, et~al., \bibinfo{year}{2022}.
\newblock \bibinfo{title}{Effect of temperature \& humdity on population dynamics of insects’ pest complex of cotton crop}.
\newblock \bibinfo{journal}{Plos one} \bibinfo{volume}{17}, \bibinfo{pages}{e0263260}.
\bibitem[{Braness et~al.(1991)Braness, Coster and Bennett}]{tab1-3}
\bibinfo{author}{Braness, G.}, \bibinfo{author}{Coster, D.}, \bibinfo{author}{Bennett, G.}, \bibinfo{year}{1991}.
\newblock \bibinfo{title}{Logistic models describing effects of temperature and humidity on residual effectiveness of chlorpyrifos and cyfluthrin formulations against {German} cockroaches ({Dictyoptera}: Blattellidae)}.
\newblock \bibinfo{journal}{Journal of economic entomology} \bibinfo{volume}{84}, \bibinfo{pages}{1746--1752}.
\bibitem[{Change(2023)}]{RN4}
\bibinfo{author}{Change, I.P.O.C.}, \bibinfo{year}{2023}.
\newblock \bibinfo{title}{Technical Summary}. \bibinfo{publisher}{Cambridge University Press}, \bibinfo{address}{Cambridge, United Kingdom and New York, NY, USA}.
\newblock pp. \bibinfo{pages}{33--144}.
\newblock \DOIprefix\doi{10.1017/9781009157896.002}.
\bibitem[{Chen et~al.(2022a)Chen, Li, Tai, Chen and Huang}]{chen2022pest}
\bibinfo{author}{Chen, C.J.}, \bibinfo{author}{Li, Y.S.}, \bibinfo{author}{Tai, C.Y.}, \bibinfo{author}{Chen, Y.C.}, \bibinfo{author}{Huang, Y.M.}, \bibinfo{year}{2022}a.
\newblock \bibinfo{title}{Pest incidence forecasting based on internet of things and long short-term memory network}.
\newblock \bibinfo{journal}{Applied Soft Computing} \bibinfo{volume}{124}, \bibinfo{pages}{108895}.
\bibitem[{Chen et~al.(2022b)Chen, Song, Yu, Zhang, Zhang, Bi, Li, Ren, Zhu, Yao et~al.}]{chen2022identifying}
\bibinfo{author}{Chen, L.}, \bibinfo{author}{Song, B.}, \bibinfo{author}{Yu, C.}, \bibinfo{author}{Zhang, J.}, \bibinfo{author}{Zhang, J.}, \bibinfo{author}{Bi, R.}, \bibinfo{author}{Li, X.}, \bibinfo{author}{Ren, X.}, \bibinfo{author}{Zhu, Y.}, \bibinfo{author}{Yao, D.}, et~al., \bibinfo{year}{2022}b.
\newblock \bibinfo{title}{Identifying soybean pod borer (\textit{{L}eguminivora glycinivorella}) resistance {QTLs} and the mechanism of induced defense using linkage mapping and {RNA}-seq analysis}.
\newblock \bibinfo{journal}{International Journal of Molecular Sciences} \bibinfo{volume}{23}, \bibinfo{pages}{10910}.
\bibitem[{Deutsch et~al.(2018)Deutsch, Tewksbury, Tigchelaar, Battisti, Merrill, Huey and Naylor}]{deutsch2018increase}
\bibinfo{author}{Deutsch, C.A.}, \bibinfo{author}{Tewksbury, J.J.}, \bibinfo{author}{Tigchelaar, M.}, \bibinfo{author}{Battisti, D.S.}, \bibinfo{author}{Merrill, S.C.}, \bibinfo{author}{Huey, R.B.}, \bibinfo{author}{Naylor, R.L.}, \bibinfo{year}{2018}.
\newblock \bibinfo{title}{Increase in crop losses to insect pests in a warming climate}.
\newblock \bibinfo{journal}{Science} \bibinfo{volume}{361}, \bibinfo{pages}{916--919}.
\bibitem[{Fan and Wang(1998)}]{tab1-1}
\bibinfo{author}{Fan, M.}, \bibinfo{author}{Wang, K.}, \bibinfo{year}{1998}.
\newblock \bibinfo{title}{Optimal harvesting policy for single population with periodic coefficients}.
\newblock \bibinfo{journal}{Mathematical biosciences} \bibinfo{volume}{152}, \bibinfo{pages}{165--178}.
\bibitem[{Fei et~al.(2024)Fei, Cui, Zhu, Xia, Xing, Gao and Shi}]{fei2024integrative}
\bibinfo{author}{Fei, H.}, \bibinfo{author}{Cui, J.}, \bibinfo{author}{Zhu, S.}, \bibinfo{author}{Xia, Y.}, \bibinfo{author}{Xing, Y.}, \bibinfo{author}{Gao, Y.}, \bibinfo{author}{Shi, S.}, \bibinfo{year}{2024}.
\newblock \bibinfo{title}{Integrative analyses of transcriptomics and metabolomics in immune response of \textit{Leguminivora glycinivorella} {Mats} to \textit{{Beauveria} bassiana} infection}.
\newblock \bibinfo{journal}{Insects} \bibinfo{volume}{15}, \bibinfo{pages}{126}.
\bibitem[{Gao et~al.(2018)Gao, Shi, Xu and Cui}]{gao2018current}
\bibinfo{author}{Gao, Y.}, \bibinfo{author}{Shi, S.}, \bibinfo{author}{Xu, M.}, \bibinfo{author}{Cui, J.}, \bibinfo{year}{2018}.
\newblock \bibinfo{title}{Current research on soybean pest management in {China}.}
\newblock \bibinfo{journal}{Oil Crop Science} .
\bibitem[{Garrido-Jurado et~al.(2011)Garrido-Jurado, Valverde-Garc{\'\i}a and Quesada-Moraga}]{tab1-4}
\bibinfo{author}{Garrido-Jurado, I.}, \bibinfo{author}{Valverde-Garc{\'\i}a, P.}, \bibinfo{author}{Quesada-Moraga, E.}, \bibinfo{year}{2011}.
\newblock \bibinfo{title}{Use of a multiple logistic regression model to determine the effects of soil moisture and temperature on the virulence of entomopathogenic fungi against pre-imaginal mediterranean fruit fly ceratitis capitata}.
\newblock \bibinfo{journal}{Biological Control} \bibinfo{volume}{59}, \bibinfo{pages}{366--372}.
\bibitem[{Hu et~al.(2012)Hu, He, Zhou, Feng and Zhang}]{hu2012synthesis}
\bibinfo{author}{Hu, D.H.}, \bibinfo{author}{He, J.}, \bibinfo{author}{Zhou, Y.W.}, \bibinfo{author}{Feng, J.T.}, \bibinfo{author}{Zhang, X.}, \bibinfo{year}{2012}.
\newblock \bibinfo{title}{Synthesis and field evaluation of the sex pheromone analogues to soybean pod borer \textit{Leguminivora glycinivorella}}.
\newblock \bibinfo{journal}{Molecules} \bibinfo{volume}{17}, \bibinfo{pages}{12140--12150}.
\bibitem[{Jiang-sheng et~al.(2021)Jiang-sheng, Jing-yi and Xia-ping}]{jiang2021hyperspectral}
\bibinfo{author}{Jiang-sheng, G.}, \bibinfo{author}{Jing-yi, F.}, \bibinfo{author}{Xia-ping, F.}, \bibinfo{year}{2021}.
\newblock \bibinfo{title}{Hyperspectral imaging for detection of \textit{Leguminivora glycinivorella} based on 3d few-shot meta-learning model}.
\newblock \bibinfo{journal}{Spectroscopy and Spectral Analysis} \bibinfo{volume}{41}, \bibinfo{pages}{2171--2174}.
\bibitem[{Kumar et~al.(2023)Kumar, Chowdary and Sarada}]{kumar2023seasonal}
\bibinfo{author}{Kumar, G.S.}, \bibinfo{author}{Chowdary, L.R.}, \bibinfo{author}{Sarada, O.}, \bibinfo{year}{2023}.
\newblock \bibinfo{title}{Seasonal incidence of major insect pests of groundnut and their natural enemies in relation to meteorological parameters.}
\newblock \bibinfo{journal}{ANGRAU} .
\bibitem[{Kuzmin et~al.(2020)Kuzmin, Anisimov and Malashonok}]{kuzmin2020soybean}
\bibinfo{author}{Kuzmin, A.}, \bibinfo{author}{Anisimov, N.}, \bibinfo{author}{Malashonok, A.}, \bibinfo{year}{2020}.
\newblock \bibinfo{title}{Soybean moth \textit{Leguminivora glycinivorella} (\textit{Lepidoptera: Tortricidae}): harmfulness in the conditions of the south of the {Amur} region}, in: \bibinfo{booktitle}{IOP Conference Series: Earth and Environmental Science}, \bibinfo{organization}{IOP Publishing}. p. \bibinfo{pages}{012026}.
\bibitem[{Le et~al.(2006)Le, Ishitani, Komai, Yamamoto and Ando}]{le2006sex}
\bibinfo{author}{Le, V.V.}, \bibinfo{author}{Ishitani, M.}, \bibinfo{author}{Komai, F.}, \bibinfo{author}{Yamamoto, M.}, \bibinfo{author}{Ando, T.}, \bibinfo{year}{2006}.
\newblock \bibinfo{title}{Sex pheromone of the soybean pod borer, \textit{Leguminivora glycinivorella} (\textit{Lepidoptera: Tortricidae}): Identification and field evaluation}.
\newblock \bibinfo{journal}{Applied entomology and zoology} \bibinfo{volume}{41}, \bibinfo{pages}{507--513}.
\bibitem[{Liu et~al.(preprint)Liu, Zheng, Li and Yu}]{tab1-2}
\bibinfo{author}{Liu, Z.}, \bibinfo{author}{Zheng, B.}, \bibinfo{author}{Li, J.}, \bibinfo{author}{Yu, J.}, \bibinfo{year}{preprint}.
\newblock \bibinfo{title}{Dynamical analysis of a pest control model with the cumulative lethal effect of periodic pesticide spraying}.
\newblock \bibinfo{journal}{Available at SSRN 5078116} .
\bibitem[{Ma and Ma(2022)}]{ma2022potential}
\bibinfo{author}{Ma, G.}, \bibinfo{author}{Ma, C.S.}, \bibinfo{year}{2022}.
\newblock \bibinfo{title}{Potential distribution of invasive crop pests under climate change: incorporating mitigation responses of insects into prediction models}.
\newblock \bibinfo{journal}{Current Opinion in Insect Science} \bibinfo{volume}{49}, \bibinfo{pages}{15--21}.
\bibitem[{Malaguit et~al.(2023)Malaguit, Mendoza, Tubay and Mata}]{malaguit2023identifying}
\bibinfo{author}{Malaguit, J.C.}, \bibinfo{author}{Mendoza, V.M.P.}, \bibinfo{author}{Tubay, J.M.}, \bibinfo{author}{Mata, M.A.E.}, \bibinfo{year}{2023}.
\newblock \bibinfo{title}{Identifying patterning behavior in a plant infestation of insect pests}.
\newblock \bibinfo{journal}{Mathematical Biosciences} \bibinfo{volume}{362}, \bibinfo{pages}{109032}.
\bibitem[{Ponti et~al.(2021)Ponti, Gutierrez, de~Campos, Desneux, Biondi and Neteler}]{ponti2021biological}
\bibinfo{author}{Ponti, L.}, \bibinfo{author}{Gutierrez, A.P.}, \bibinfo{author}{de~Campos, M.R.}, \bibinfo{author}{Desneux, N.}, \bibinfo{author}{Biondi, A.}, \bibinfo{author}{Neteler, M.}, \bibinfo{year}{2021}.
\newblock \bibinfo{title}{Biological invasion risk assessment of tuta absoluta: mechanistic versus correlative methods}.
\newblock \bibinfo{journal}{Biological Invasions} \bibinfo{volume}{23}, \bibinfo{pages}{3809--3829}.
\bibitem[{Shi et~al.(2018)Shi, Cui, Zhu, Xu and Wang}]{shishusen}
\bibinfo{author}{Shi, S.}, \bibinfo{author}{Cui, J.}, \bibinfo{author}{Zhu, S.}, \bibinfo{author}{Xu, W.}, \bibinfo{author}{Wang, X.}, \bibinfo{year}{2018}.
\newblock \bibinfo{title}{Genetic differentiation of soybean pod borer in different geographical populations based on mitochondrial {COI} gene sequence}.
\newblock \bibinfo{journal}{Acta Phytophylacica Sinica} \bibinfo{volume}{45}, \bibinfo{pages}{214--222}.
\bibitem[{Skend{\v{z}}i{\'c} et~al.(2021)Skend{\v{z}}i{\'c}, Zovko, {\v{Z}}ivkovi{\'c}, Le{\v{s}}i{\'c} and Lemi{\'c}}]{skendvzic2021impact}
\bibinfo{author}{Skend{\v{z}}i{\'c}, S.}, \bibinfo{author}{Zovko, M.}, \bibinfo{author}{{\v{Z}}ivkovi{\'c}, I.P.}, \bibinfo{author}{Le{\v{s}}i{\'c}, V.}, \bibinfo{author}{Lemi{\'c}, D.}, \bibinfo{year}{2021}.
\newblock \bibinfo{title}{The impact of climate change on agricultural insect pests}.
\newblock \bibinfo{journal}{Insects} \bibinfo{volume}{12}, \bibinfo{pages}{440}.
\bibitem[{Sokame et~al.(2021)Sokame, Tonnang, Subramanian, Bruce, Dubois, Ekesi and Calatayud}]{sokame2021system}
\bibinfo{author}{Sokame, B.M.}, \bibinfo{author}{Tonnang, H.E.}, \bibinfo{author}{Subramanian, S.}, \bibinfo{author}{Bruce, A.Y.}, \bibinfo{author}{Dubois, T.}, \bibinfo{author}{Ekesi, S.}, \bibinfo{author}{Calatayud, P.A.}, \bibinfo{year}{2021}.
\newblock \bibinfo{title}{A system dynamics model for pests and natural enemies interactions}.
\newblock \bibinfo{journal}{Scientific reports} \bibinfo{volume}{11}, \bibinfo{pages}{1401}.
\bibitem[{Sun et~al.(2018)Sun, Pu and Lin}]{tab1-6}
\bibinfo{author}{Sun, S.}, \bibinfo{author}{Pu, L.}, \bibinfo{author}{Lin, Z.}, \bibinfo{year}{2018}.
\newblock \bibinfo{title}{Dynamics of the logistic harvesting model with infinite delay on periodically evolving domains}.
\newblock \bibinfo{journal}{Commun. Math. Biol. Neurosci.} \bibinfo{volume}{2018}, \bibinfo{pages}{Article--ID}.
\bibitem[{Tang and Cheke(2008)}]{tab1-5}
\bibinfo{author}{Tang, S.}, \bibinfo{author}{Cheke, R.A.}, \bibinfo{year}{2008}.
\newblock \bibinfo{title}{Models for integrated pest control and their biological implications}.
\newblock \bibinfo{journal}{Mathematical Biosciences} \bibinfo{volume}{215}, \bibinfo{pages}{115--125}.
\bibitem[{Wakamura(1992)}]{wakamura1992development}
\bibinfo{author}{Wakamura, S.}, \bibinfo{year}{1992}.
\newblock \bibinfo{title}{Development in application of synthetic sex pheromone to pest management}.
\newblock \bibinfo{journal}{Jpn. Pestic. Inform.} \bibinfo{volume}{61}, \bibinfo{pages}{26--31}.
\bibitem[{Yang et~al.(2024)Yang, Wang, Ding, Li and Zhang}]{yang2024predicting}
\bibinfo{author}{Yang, M.}, \bibinfo{author}{Wang, Y.}, \bibinfo{author}{Ding, W.}, \bibinfo{author}{Li, H.}, \bibinfo{author}{Zhang, A.}, \bibinfo{year}{2024}.
\newblock \bibinfo{title}{Predicting habitat suitability for the soybean pod borer \textit{Leguminivora glycinivorella} ({Matsumura}) using optimized maxent models with multiple variables}.
\newblock \bibinfo{journal}{Journal of Economic Entomology} \bibinfo{volume}{117}, \bibinfo{pages}{1796--1808}.

\end{thebibliography}
\end{document}